\documentclass[MSNbibl,number,citesort,dvips]{arxbj}
\usepackage{graphicx}

\aid{0}
\volume{18}
\issue{3}
\pubyear{2012}
\firstpage{1031}
\lastpage{1041}
\doi{10.3150/10-BEJ350}

\makeatletter
\newcommand{\prob}{\mathbf{P}}
\newcommand{\mean}{\mathbf{E}}
\newcommand{\ind}{\mathbf{1}}
\newcommand{\bern}{\operatorname{{\textsf{Bern}}}}
\newcommand{\pois}{\operatorname{{\textsf{Pois}}}}
\newcommand{\der}{\mathrm{d}}
\newcommand{\kernel}{\mathbf{K}}
\newtheorem{theorem}{Theorem}[section]
\newcommand{\iid}{\mathop{\sim}\limits^{\mathrm{i.i.d.}}}
\newcommand{\pswap}{{p_{\mathrm{swap}}}}
\makeatother

\begin{document}
\begin{frontmatter}

\title{Spatial birth--death swap chains}
\runtitle{Spatial birth--death swap chains}

\begin{aug}
\author{\fnms{Mark} \snm{Huber}\corref{}\ead[label=e1]{mhuber@cmc.edu}}
\runauthor{M. Huber}
\address{Department of Mathematics and Computer Science, 850 Columbia Ave.,
Claremont McKenna College, Claremont, CA  91711, USA. \printead{e1}}
\end{aug}

\received{\smonth{6} \syear{2008}}
\revised{\smonth{11} \syear{2010}}

%
\begin{abstract}
Markov chains have long been used for generating random variates from
spatial point processes. Broadly
speaking, these chains fall into two categories: Metropolis--Hastings
type chains running in discrete time and spatial
birth--death chains running in continuous time.
These birth--death chains only allow
for removal of a point or addition of a point. In this paper
it is shown that the addition of transitions where a point is moved from
one location to the other can aid in shortening the mixing time of the
chain. Here the mixing time of the chain is analyzed
through coupling, and use of the swap moves allows for analysis of a
broader class of chains. Furthermore, these swap moves can be
employed in perfect sampling algorithms via the dominated coupling from
the past procedure of Kendall and M{\o}ller.
This method can be applied to any
pairwise interaction model with repulsion. In particular,
an application
to the Strauss process is developed in detail, and the swap
chains are shown
to be much faster than standard birth--death chains.
\end{abstract}

%

\end{frontmatter}

\section{Introduction}

Spatial point processes are in wide use in statistical modeling
(see \cite{illianpss2008} for an overview). Typically
finite point processes are modeled as being
absolutely continuous with respect to a Poisson
point process.
That is, they have a density $f(x) / c$ where
$f(x)$ is an easily computable function but the normalizing constant
$c$ of the density is impractical to compute.
A Monte Carlo algorithm gains information about $f(x) / c$
by studying random variates drawn from the distribution
the density describes.

To obtain these variates, a Markov chain is built
whose stationary distribution matches the target distribution.
Metropolis--Hastings chains run in discrete time
(see \cite{geyer1999}), and
the spatial birth--death chain approach of Preston \cite{preston1977}
runs in continuous time. In \cite{cliffordn1994} problems were given
where the Metropolis--Hastings approach is faster
than Preston's.

The drawback of
these Markov chain Monte Carlo methods is that
unless the mixing time of the Markov chain is known,
the quality of the variates is suspect. Heuristics such
as the autocorrelation test can
prove that a chain has not mixed, but cannot establish
the positive claim that a chain has mixed.

Perfect simulation algorithms solve this problem. They
generate samples exactly from the desired distribution without the need
to know the mixing time of a Markov chain.
Kendall \cite{kendall1995} showed how the
coupling from the past (CFTP) idea of Propp and Wilson~\cite{proppw1996}
could be used together with a spatial birth and death chain to obtain
samples from area interaction processes. Kendall and
M{\o}ller \cite{kendallm2000} showed how this method could be extended
to any locally stable point process using a method they called
\textit{dominated CFTP}.
They also considered perfect sampling using
Metropolis--Hastings chains,
but restricted these chains to only adding or deleting a point at each step.

So \cite{cliffordn1994} indicates that
Metropolis--Hastings chains can beat continuous time chains, but~\cite{kendallm2000} shows how
to exactly sample using continuous time chains. The
goal of this work is to introduce a new swap move to the continuous time
chains that speeds up convergence, while still allowing
for perfect simulation.

In Section~\ref{SEC:swap_move} the
theory behind spatial birth--death chains with the new swap move is developed,
and an example of such a chain is given for the Strauss process.
Section~\ref{SEC:dcftp} reviews the use of dominated coupling
from the past, and shows how the addition
of swap moves fits into this protocol.
Section~\ref{SEC:runtime} bounds the
expected running time of the procedure for a restricted class of models.

\section{Spatial point processes}
\label{SEC:swap_move}
Dyer and Greenhill \cite{dyerg2000b} first introduced a swap move
for hard core point processes in discrete spaces. In this section their method
is extended to more general
point processes.

For ease of exposition, we consider here point processes that do not
contain multiple points. Let $S$ be a separable measurable set,
and $\lambda$ be a diffuse measure on $S$ (so
$\lambda(\{v\}) = 0$ for all $v \in S$) such that $\lambda(S) < \infty$.
(Typically $S$ is a bounded Borel set of $\mathbf{R}^2$.)
Then a~Poisson point process
is a finite subset of $S$ chosen as follows. First, let $N$ be a Poisson
distributed random variable with parameter $\lambda(S)$. Next, let
$X_1,\ldots,X_N$ be independently and identically distributed (i.i.d.)
and drawn from the probability distribution $\lambda(\cdot)/\lambda(S)$.
Then $\{X_1,\ldots,X_N\}$ (called a \textit{configuration}) is a draw
from a Poisson point process with intensity measure $\lambda(\cdot)$
over $S$.
Let $\mu$ be the distribution of the configuration and $\Omega$ the
set of all possible configurations. More details of $\mu$
and $\Omega$ can be found in \cite{carterp1972,preston1977}.

As an example of data modeled using these types of processes,
Harkness and Isham~\cite{harknessi1983}
studied locations of ant nests in a rectangular region $R$.
With two types of ants,
$S = R \times\{0,1\}$ and $\lambda$ is the product
of Lebesgue and a measure on $\{0,1\}$.

The processes considered here are absolutely continuous with respect to
$\mu$ with density~$f$ satisfying a \textit{local stability condition}
(as in \cite{kendallm2000}):
%
%
\begin{equation}
\label{EQN:local_stability}
(\exists K > 0)(\forall x \in\Omega)(\forall v \in S \setminus x)
\bigl(f(x \cup\{v\}) \leq K f(x)\bigr).
\end{equation}
Many point processes of interest meet this condition,
including the area interaction
process \cite{baddeleyvl1995,widomr1970}, the
Strauss process \cite{kellyr1976,strauss1975}
and the continuous random cluster
model \cite{haggstromvlm1999}.

\subsection{Spatial birth--death swap chains}
\label{SEC:spatial_chains}

The development of the swap move given here follows the framework of
Preston \cite{preston1977}, who introduced the use of
spatial birth--death chains for these problems. These chains are examples
of \textit{jump processes}, where at a given state $x$, the chain stays
in the state for an exponential length of time with expected value given
by $1 / \alpha(x)$. The state then jumps to a new state using
kernel~$\kernel$, so the probability that the new state is
in $A$ is $\kernel(x,A)$), independent of the past history
(see \cite{feller1966}, Chapter X, for the details of jump processes).

In the Preston framework,
the rate of births (addition of points to the configuration) and
deaths (deletion of points from the configuration) depends only on
the current state:
\begin{itemize}
\item There exists a non-negative measurable birth rate function $b$
from $\Omega\times S$ equipped with the standard
product $\sigma$-field to $\mathbf R$ with the Borel $\sigma$-field.
Call $b(x,v)$ the \textit{birth rate} at which point $v$ is added to
configuration $x$.
\item There exists a
non-negative measurable death rate function
$d$ from $\Omega\times S$ equipped with the standard product $\sigma$-field
to $\mathbf{R}$ with the Borel $\sigma$-field. Furthermore,
$w \in x \Rightarrow d(x,w) > 0$ and
$w \notin x \Rightarrow d(x,w) = 0$. Then $d(x,w)$ is the \textit{death rate}
at which a~point $w$ is removed from configuration $x$.
\end{itemize}
To this birth--death framework we now add a swap rate:
\begin{itemize}
\item There exists a non-negative measurable swap rate function
$s$ from $\Omega\times S \times S$
equipped with the standard product $\sigma$-field
to $\mathbf{R}$ with the Borel $\sigma$-field. Furthermore,
$w \notin x \Rightarrow s(x,w,v) = 0$. So $s(x,w,v)$ is the \textit
{swap rate}
at which point $w$ is removed and point $v$ is added.
\end{itemize}
The birth, death, and swap rates are used to build a kernel $\kernel$ for
the Markov chain as follows.
For all $A \in\mathcal{B}$, let
$K_{b}(x,A) = \int_{v \in S} b(x,v) \ind(x \cup\{v\} \in A)  \lambda
(\der v).$
When $K_b(x,\Omega) < \infty$ for all $x$ in $\Omega$, the birth kernel
is
$\kernel_b(x,A) = K_b(x,A) / K_b(x,\Omega).$
Similarly, $K_d(x,A) = \sum_{w \in x} d(x,w) \ind(x \setminus\{w\} \in A),$
which always has
a finite number of terms and so
$\kernel_d(x,A) = K_d(x,A) / K_d(x,\Omega).$ The total rate of births is
$r_b(x) = \int_{v \in S} b(x,v)  \lambda(\der v)$, and the total rate of
deaths is $r_d(x) = \sum_{v \in x} d(x,v)$.

For the swap kernel, set
$K_s(x,A) = \sum_{w \in x}
\int_{v \in S} s(x,w,v) \ind(x \cup\{v\} \setminus\{w\} \in A) \lambda(\der v).$
When $K_s(x,\Omega) < \infty$ for all $x \in\Omega$, let
%
%
\begin{equation}
\kernel_s(x,A) = K_s(x,A) / K_s(x,\Omega),\qquad  r_s(x) = \sum_{w \in x} \int_{v \in S} s(x,w,v)  \lambda(\der v).
\end{equation}
The overall rate at which the configuration changes is
$\alpha(x) = r_b(x) + r_d(x) + r_s(x),$
and the overall kernel is:
%
%
\begin{equation}
\kernel(x,A) = \kernel_b(x,A) \frac{r_b(x)}{\alpha(x)} +
\kernel_d(x,A) \frac{r_d(x)}{\alpha(x)} +
\kernel_s(x,A) \frac{r_s(x)}{\alpha(x)}.
\end{equation}

Harris recurrence guarantees that a Markov process has a unique
invariant measure (see \cite{azemakr1967} for
details of Harris recurrence
in the continuous-time context).
Kaspi and Mandelbaum \cite{kaspim1994}
showed that a continuous-time chain is Harris recurrent if and only if
there exists a non-zero $\sigma$-finite measure where $X$ almost surely
hits sets with positive measure.

In particular, for all the chains considered here, the death rate
equals the number of points in the configuration, and the birth rate
is bounded above by a constant. This forces the chain to visit the
empty configuration infinitely often, making it Harris recurrent.

The detailed balance conditions (that imply $f$ is invariant)
for jump processes are:
$f(x) \alpha(x) \times \kernel(x,\der y)  \,\der\mu(x) =
f(y) \alpha(y) \kernel(y,\der x)  \,\der\mu(y).$
For moves from configurations with~$n$ points to those with
$n + 1$ (or vice versa),
the detailed balance conditions are
satisfied \cite{preston1977,ripley1977} when
the rate of births balance the rate of deaths with respect to
$f$.
So
%
%
\begin{equation}
\label{EQN:birth_death_rate}
f(x) b(x,v) = f(x \cup\{v\})d(x \cup\{v\},v).
\end{equation}

Swap moves stay inside the same dimensional space, and it is
straightforward to show that reversibility for swap moves holds when
%
%
\begin{equation}
\label{EQN:swap_rate}
f(x)s(x,w,v) = f(x \cup\{v\} \setminus\{w\})s(x \cup\{v\} \setminus\{
w\},v,w).
\end{equation}

\subsection{Locally stable repulsive point processes}

Kendall and M{\o}ller \cite{kendallm2000} describe how to create a
jump process with stationary density $f$ for locally stable processes.
Briefly, their method works
as follows. Two coupled chains will be run: the \textit{dominating chain}
with state $D(t)$ at time $t$ and the \textit{target chain} with state
$X(t)$ at
time $t$. It will always be true that $X(t) \subseteq D(t)$.
Each point $w \in D(t)$ has death rate $d(D(t),w) = 1$. If a point dies
that is also in~$X(t)$, it is removed from both~$X(t)$ and $D(t)$.
The rate of births for the dominating chain is $r_b = K \lambda(S)$, where
$K$ is the local stability
constant in equation~(\ref{EQN:local_stability}). If a birth
occurs, a point $v$ is chosen according to the probability measure
$\lambda(\cdot) / \lambda(S)$. Then~$v$ is always added to
$D(t)$ to get the next dominating state, but is only added to
$X(t)$ with probability $f(X(t) \cup\{v\}) / [K f(X(t))].$
Assume that each point
$v$ born in $D(t)$ is marked with a uniform draw from $[0,1]$. Then the
point is born in $X(t)$ if the mark falls below
$f(X(t) \cup\{v\}) / [K f(X(t))].$

Suppose $X(0) \subseteq D(0)$. Then
since deaths are always accepted in both chains, but a birth in the
dominating chain might not occur in the target chain, the dominating
configuration will be a superset of the target configuration for all $t
\geq0$.

Adding a swap move to this birth death framework can be done
automatically
when the rejection of a birth $v$ can be linked to a single point
$w \in X(t)$. Consider an example.

\subsubsection*{Strauss model}
In the Strauss model
\cite{kellyr1976,strauss1975}, the density has a factor that is
exponential in the number of pairs of points that lie within distance
$R$ of each other. Let $\rho$ be a metric on $S$ (usually Euclidean
distance), then the density can be written:
%
%
\begin{equation}
f_{S}(x) = Z^{-1}_{(\beta_1,\beta_2,R)} \beta_1^{\# x} \beta_2^{s(x)},
\qquad  s(x) = \sum_{\{v,v'\}:v \in x, v' \in x \setminus\{v\}} \ind\bigl(\rho
(v,v') \leq R\bigr),
\end{equation}
where $Z_{(\beta_1,\beta_2,R)}$ is the normalizing constant for the
density. As noted in \cite{kellyr1976},
in order for $Z_{(\beta_1,\beta_2,R)}$ to be finite (and
hence for the density to exist) $\beta_2$ must be at most 1.
In addition, \cite{kellyr1976} generalizes the Strauss process to the
pairwise interaction process. All methods presented here are written
for the Strauss process for simplicity, but work equally well for the
pairwise interaction process.

Let $x$ be the state of the target chain, and suppose point
$v$ is born in the dominating chain. Call point $w \in x$
a neighbor of $v$ if $\rho(v,w) \leq R$. The Strauss process
is locally stable with $K = \beta_1$, so
the chance of accepting $v$ into $x$ is
$f(x \cup\{v\}) / [K f(x)] = \beta_2^{s(x,v)},$ where
$s(x,v) = \sum_{v' \in x} \ind(\rho(v',v) \leq R)$ is the number of
neighbors of $v$ in $x$.

Let $\bern(p)$ denote the Bernoulli distribution with parameter $p$.
One way to draw $B \sim\bern(\beta_2^{s(x,v)})$ is to draw
$B_1,\ldots,B_{s(x,v)} \iid\bern(\beta_2)$ and
set $B = B_1B_2\cdots B_{s(x,v)}.$
(Here $\iid$ denotes that the
draws are to be independent and identically distributed.)

When $B_i = 0$, say that the point indexed by $i$ \textit{blocks} the birth
of $v$.
Suppose that $v$ is blocked by a single neighbor $w$.
Then the swap move removes $w$, and allows the
birth of $v$.
Call this new configuration $x'$.
The probability of swapping from $x$ to $x'$ (given birth $v$) is
$\beta_2^{s(x,v) - 1}(1 - \beta_2)$. This makes it straightforward
to check
that $f(x)s(x,w,v) = f(x')s(x',v,w)$, so
(\ref{EQN:swap_rate}) is satisfied.
To implement this swap move, simply mark each point~$v$
born in
$D(t)$ with an i.i.d. sequence of $\bern(\beta_2)$ random
variables.

\section{Perfect simulation by dominated CFTP}
\label{SEC:dcftp}

In the previous section it was shown how to couple a dominating
chain and target chain using standard birth--death chains and the
new birth--death swap chain. Here a further coupling is built that
allows exact draws to be taken from the
stationary distribution of the target chain using the
dominating CFTP (dCFTP)
method of Kendall
and M{\o}ller~\cite{kendallm2000}.

Both $X(t)$ and $D(t)$ are time-reversible, so
they can be run backwards in time as easily as forwards while
maintaining the property that if $X(0) \subseteq D(0)$, then
$X(t) \subseteq D(t)$ for all $t \in(-\infty,0]$.
(A more detailed introduction to dCFTP can be found
in \cite{kendallm2000}.)

So far two chains (the dominating and target) have been coupled,
but now consider two more chains, called the \textit{lower chain}
and \textit{upper chain}, denoted $L(t)$ and $U(t),$ respectively.
Suppose that these four chains have the \textit{sandwiching property} that
%
%
\begin{equation}
\label{EQN:bounding}
L(t) \subseteq X(t) \subseteq U(t) \subseteq D(t)\qquad
\mbox{for all }
t \in(-\infty,0].
\end{equation}
The process $(L(t),U(t))$ can also be thought of as a bounding
process for $X(t)$ (see \cite{huber2004a}).

Suppose $X(0)$ is drawn from the stationary distribution.
Then if $L(0) = U(0)$, $X(0)$ also equals the lower and upper
chain, and the state they all equal is a draw from the stationary
distribution. This is the idea behind CFTP.

For each positive integer $N$, a lower and upper chain can be created.
Consider $D(t)$ moving backward through time, and let $\tau_N$ denote
the time where the $N$th backward event occurs. Set $L_N(\tau_N)$ to
the empty
configuration, and $U_N(\tau_N) = D(\tau_N)$.

Every time there is an event at time $t$
(either a birth or death in the
dominating process moving forwards in time) it is important to ensure
that $U_N(t)$ and $L_N(t)$ continue to
bound $X(t)$ once the event updates the chain. That is, if a point
$v$ is added to the target chain state, it must also be added to
the upper chain. If a point $w$ is removed from the target chain
state, it must also be removed from the lower chain.
Such a coupling has
the \textit{funneling property} (see \cite{berthelsenm2002}). All
the couplings used here have this important property.

An induction
argument shows that the funneling property implies
$L_N(0) \subseteq X(0) \subseteq U_N(0)$. Note
if $L_N(0) = U_N(0)$, then $X(0)$ is trapped between them and also
equals this common value. This is the coupling part of CFTP.

The ``from the past'' part of CFTP works as follows. Suppose
$L_N(0) \not= U_N(0)$. Then increase the value of $N$ and try
again. Let $N' > N$.
The first $N$ events for the dominating process (looking backward
in time from time 0)
have already been generated, these same events must be used in
subsequent evaluations of the bounding process. Therefore, only
$N' - N$ additional events need to be generated. Once these
events have been generated, run $L_{N'}$ and $U_{N'}$ forward until
$L_{N'}(0)$ and $U_{N'}(0)$ can be compared.

If
$L_{N}(0) = U_N(0)$ for some $N$, then
$L_{N'}(0) = U_{N'}(0)$ for all $N' > N$ as well, so
it is not necessary to try every value of $N$. Propp and
Wilson \cite{proppw1996} noted that by
doubling $N$ at each step, the total number of checked events
is at most twice the minimum number.
The choice of $N_{\mathrm{initial}}$ is arbitrary, but $L_N(0)$
cannot equal $U_N(0)$ unless every point in $D(\tau_N)$ has died
by time 0. For simplicity, here $N_{\mathrm{initial}}$ is set equal to the expected
number of points in the dominating process at time 0, which is $K
\lambda(S)$
(see \cite{berthelsenm2002} for a more advanced approach to choosing
$N_{\mathrm{initial}}$).

Kendall and M{\o}ller showed (Theorem 2.1 of \cite{kendallm2000})
that as long as the probability
that~$D(t)$ visits the empty configuration in $[0,t]$ goes to 1 as
$t$ goes to infinity, this procedure will terminate in finite time
with probability 1. The resulting configuration $L_N(0) = U_N(0)$
is a draw exactly from the target distribution.

Now consider the question: How should
the lower and upper chains be updated for each event in the
dominating process so the funneling property holds for the swap move?

\subsection{Updating the bounding process}
\label{SBS:update}

For a jump process $A(t)$, let $A(t-)$ denote the limit as $\varepsilon$
goes to 0 of $A(t - \varepsilon)$, that is, the state of the process right
before time $t$. The bounding process needs to be updated if a point
is born or dies at time $t$. The procedure followed is the same as
given in \cite{huber2004a}.

If a point $w \in X(t-)$ dies, it is removed from $X(t)$, and so can be
removed from both~$L_N(t)$ and $U_N(t)$. Now suppose point $v$ is born
into the dominating chain at time~$t$.

Case 1:
Point $v$ is blocked by at most one point $w$ in $U_N(t-)$. Then
$X(t-) \subseteq U_N(t-)$ and so if $w \in X(t)$, then $w$ is swapped
away by $v$,
and if $w \notin X(t-)$, then $v$ can be born. So either way
$X(t) = X(t-) \setminus\{w\} \cup\{v\}$, $w$ is removed from $U_N(t)$
(and $L_N(t)$ if it is there also)
and $v$ is added to both $L_N(t)$ and $U_N(t)$.

Case 2: The point $v$ is blocked by at least
two points in $L_N(t-)$. Then
there are at least two blocking points
in $X(t-)$, so the birth does
not occur in $L_N(t), X(t)$ or $U_N(t)$.

Case 3: the point $v$ is blocked by at most one point in $L_N(t-)$,
and at least two points in $U_N(t-)$. Then if $X(t-)$ contains the
two blocking points in $U_N(t-)$, the swap does not occur, but
if it only contains the single blocking point in $L_N(t-)$, the
swap does occur. The result is that the birth $v$ must be added
to $U_N(t)$ (but not to $L_N(t)$)
to ensure $X(t) \subseteq U_N(t)$, and any
blocking point in $L_N(t-)$ must be removed from $L_N(t)$.

\begin{figure}

\includegraphics{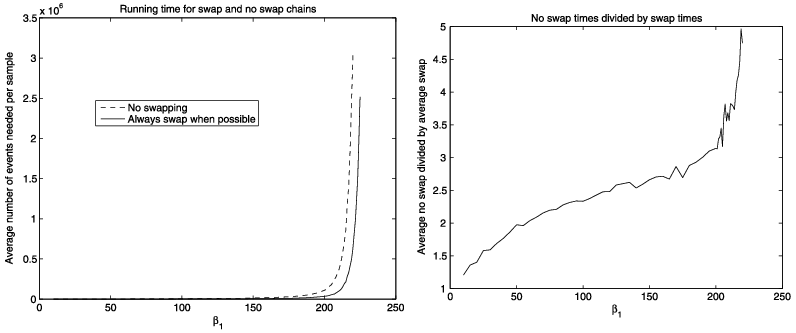}

\caption{Running time of dCFTP for Strauss model on $S = [0,1]^2$,
$\beta_2 = 0.5$, $R = 0.05$, $\lambda$ is Lebesgue measure.}
\label{FIG:hctime}
\end{figure}

Figure \ref{FIG:hctime} shows the running time advantage gained by
using the
swap move. The times
are measured in number of events generated
by dominated CFTP (dCFTP).
On the left are the raw number of times for the chain without the swap move
and with the swap move. The plot on the right
shows the ratio of these two times. Note that as $\beta_1$ gets larger,
the speedup gained by using the swap move also increases.

%

\section{Analyzing the running time}
\label{SEC:runtime}

Consider how many events must be generated before the dominated
coupling from the past procedure terminates, that is, before
$U_N(0) = L_N(0)$.
Deaths in $U_N(t) \setminus L_N(t)$
cause the bounding process to move together, while births can add
a point to $U_N(t)$ but not to~$L_N(t)$, and the swap move sometimes
removes a
point from $L_N(t)$ but not $U_N(t)$. Therefore, it is reasonable that
the perfect simulation algorithm will run faster in situations where
the birth rate is low.

In this section it is shown that, for perfect simulation of the Strauss
process, the original no swap chain takes
(with high probability) a small number of steps per perfect sample
when $\beta_1$ and $R$ are not too large, and $\beta_2$ is
not too small. By creating a mixture
of the swap chain and no swap chain, it is possible to improve
this result to where it applies for values of $\beta_1$ that
are twice as large as for the no swap chain.

The mixture works as follows: At each step,
with probability $\pswap$, the swap move chain is
used, while with probability $1 - \pswap$, the original no swap chain
is used. The best theoretical bound is achieved when $\pswap= 1/4$.

\begin{theorem}
\label{THM:runtime}
Suppose that $N$ events are generated backwards in time and then run
forward to get $U_N(0)$ and $L_N(0)$.
Let $B(v,R)$ be the area
within distance $R$ of $v \in S$, let
$r = \sup_{v \in S} \lambda(B(v,R)).$

If $\beta_1 (1 - \beta_2) r < 1$, then for the chain without the swap move
%
%
\begin{equation}
\label{EQN:runtime}
\hspace*{-3pt}\prob\bigl(U_N(0)\! \neq \!L_N(0)\bigr) \!\leq\! 2\exp(-0.09N)
\!+\! \beta_1 \lambda(S)
\exp\bigl(-N \bigl(1\! -\! \beta_1 (1 \!-\! \beta_2) r\bigr) / (4 \beta_1 \lambda(S))\bigr).
\end{equation}

If $\beta_1 (1 - \beta_2) r < 2$, then
for the chain where a swap is executed with probability $1/4$,
%
%
\begin{eqnarray}
\label{EQN:secondTheorem}
\prob\bigl(U_N(0) \neq L_N(0)\bigr) &\leq& 2\exp(-0.09N)
\nonumber
\\[-8pt]
\\[-8pt]
\nonumber
&&{}+ \beta_1 \lambda(S)
\exp\bigl(-N \bigl(1 - 0.5 \beta_1 (1 - \beta_2) r\bigr) / (4 \beta_1 \lambda(S))\bigr).
\end{eqnarray}
\end{theorem}

Why the value of 1$/$4 for the probability? This is an artifact of the proof
technique.
The theorem only gives sufficient, not necessary, conditions for the
algorithm to be fast, and simulation experiments indicate
that the algorithm actually takes the fewest steps when the swap
moves are used as often as possible (reasons why this could be true are noted
below in the proof of the theorem).

Theorem \ref{THM:runtime} has immediate consequences for the expected
running time of dominated CFTP.
Recall that in dCFTP the number of events was doubled each
time. Say $\prob(U_N(0) \neq L_N(0)) \leq a\exp(-bN),$ and
let $T$ be the number of events generated in a call of dCFTP. Then
for $T \geq t$,
dCFTP must have failed on a run of length at least $t / 2$. So
%
%
\begin{equation}
\mean[T] = \sum_{N=1}^\infty\prob(T \geq N) \leq
\Biggl[\sum_{N=1}^{\lceil(2/b)\ln a \rceil} 1 \Biggr]+
\sum_{N=\lceil(2/b)\ln a \rceil}^\infty a
\exp(-b N / 2),
\end{equation}
which makes $\mean[T] = \mathrm{O}(\ln a / b)$, and
the mean running time
$\mathrm{O}(\beta_1 \lambda(S) (\ln\beta_1 \lambda(S)))$ for the no swap chain when
$\beta_1 (1 - \beta_2) r < 1$ and in the $1/4$-swap chain when
$\beta_1 (1 - \beta_2) r < 2$.

\begin{pf*}{Proof of Theorem \protect\ref{THM:runtime}}
Recall $U_N(\tau_N) = D(\tau_N)$, a
Poisson spatial point process with parameter
$\beta_1 \lambda(S)$. $L_N(\tau_N)$ is the empty configuration, and
the bounding processes are run forward in time. Let $Q(t) = U_N(t)
\setminus L_N(t)$.
Then the chains have come together if and
only if $\#Q(0) = 0$. Begin by considering the no swap chain.

\textit{Strauss no swap move}. All individual death rates
are 1, so the total rate of deaths of points in $Q(t)$ is just $\#Q(t)$.
Call a death a \textit{good event} since it reduces $\#Q(t)$ by 1.

For $\#Q(t)$ to increase by 1 (call this a \textit{bad event}), a birth
must occur
at $v$ and be added to $U_N(t)$ but not $L_N(t)$. Let $w$ be any point
in $Q(t)$.
Then for $Q(t)$ to give rise to another point in $Q(t)$, a point $v$
must be
born within distance $R$ of $w$ and the $\bern(\beta_2)$ draw must be 0.
The area surrounding $w$ is at most $r$, and the Bernoulli draw acts as
a~thinning procedure in a Poisson process
(see Appendix G of \cite{mollerw2004}.)
So the rate at which~$w$
creates new points in $Q(t)$ is at most $\beta_1 (1 - \beta_2) r$, and
the overall rate
of bad events is at most $\beta_1 (1 - \beta_2) r \#Q(t)$.

Suppose the rate of bad events is smaller than the rate of good events.
The probability
that one event occurs in the time interval from $t$ to $t + h$ is
proportional to $h$, the probability that $n$ events occurs is
$O(h^n)$. Hence
\[
\mean\bigl[\mean[\#Q(t + h)|U(t),L(t)] - \#Q(t)\bigr] \leq
\mean\Biggl[\bigl(\# Q(t) \beta_1 (1 - \beta_2) r - \#Q(t)\bigr)h + \sum_{i=2}^\infty i\mathrm{O}(h^i)\Biggr],
\]
which means
\[
\lim_{h \rightarrow0} \frac{\mean[\mean[\#Q(t + h)|U(t),L(t)] - \#
Q(t)]}{h} \leq
-\mean\bigl[\# Q(t)\bigl(1 - \beta_1 (1 - \beta_2) r\bigr)\bigr].
\]
Let $q(t) = \mean[\#Q(t)]$, and let $\tau_N$ be the time of the $N$th
event moving backwards
in time. Then $q(\tau_N) \leq\mean[\#D(\tau_N)] = \beta_1 \lambda(S)$,
so together
with $q'(t) \leq- q(t)(1 - \beta_1 (1 - \beta_2) r)$:
\[
q(t) \leq\beta_1 \lambda(S) \exp\bigl(-t\bigl(1 - \beta_1 (1 - \beta_2) r\bigr)\bigr).
\]
By Markov's inequality, $\prob(Q(0) \neq\varnothing) = \prob(\#Q(0) \geq
1) \leq q(0)$.

Now fix $N$, the number of events to run back in time, and set $t = N /
[4 \beta_1 \lambda(S)]$.
The chance $Q(0)$ does not equal 0 starting at $-t$ is at most
$\exp(- N / [4 \beta_1 \lambda(S)](1 - \beta_1 (1 - \beta_2) r)).$

Using Chernoff bounds \cite{chernoff1952}, it can be shown that for
$A \sim\pois(\alpha),$ $\prob(A > 2 \alpha) \leq \exp(-\alpha(2 \ln2
- 2 + 1))$.
So after $t$ time, the probability that more than $N / 2$ events were
generated in a Poisson process with rate $\beta_1 \lambda(S)$ is at most
$\exp(-(N / 4)(2 \ln2 - 2 + 1)).$ Both the times of the births and
times of deaths (viewed
individually) are Poisson processes with rate $\beta_1 \lambda(S)$,
therefore the
probability that either uses more than $N / 2$ events (by the union
bound) is at most
$2 \exp(-0.09N)$. But if at this time each process used at most $N / 2$
events, then moving back
in time $N$ events puts the user even farther back in time, and if
coalescence occurs at $-t$,
it will also occur starting at $\tau_N$. Again using the union bound, the
probability of failure is at most
\[
2\exp(-0.09N) + \exp\bigl(-N\bigl(1 - \beta_1 (1 - \beta_2) r\bigr)\bigr).
\]

\textit{Strauss with swap move}. Now consider what happens when
$\pswap> 0$. The rate of good events (deaths) remains unchanged, but
the rate
of bad events changes. In Section \ref{SBS:update}, Case 1 leaves $\#Q(t)$
unchanged or reduces it by 1, Case 2 leaves $\#Q(t)$ unchanged, and
Case 3 increases $\#Q(t)$ by 1 or 2. To be precise, let $A_L$ be the
set of blocking
points in $L_N(t-)$, and $A_U$ be the set of blocking points in
$U_N(t-)$. Then
the situations that change $\#Q(t)$ are:\vspace*{9pt}

{\centering{
\begin{tabular*}{\textwidth}{@{\extracolsep{\fill}}lllll@{}}
\hline
Type &
$\#A_U$ & $\#A_L$ & $\#Q(t) - \#Q(t-)$ no swap & $\#Q(t) - \#Q(t-)$
with swap \\
\hline
1 & 1 & 0 & 1 & $-1$ \\
2 & at least 2 & 1 & 0 & \phantom{$-$}2 \\
3 & at least 2 & 0 & 1 & \phantom{$-$}1\\
\hline
\end{tabular*}}
}\vspace*{6pt}

Let $b_1$ denote the area
of the region where a birth is Type 1, with $b_2$ and $b_3$
defined similarly.
Together, the rate of change from births is:
\[
b_1[(1 - \pswap) - \pswap] + b_2[2\pswap] + b_3[(1 - \pswap) + \pswap].
\]
Any point in $b_3$ neighbors at least two points in $\#Q(t-)$, and
points in $b_1$ or $b_2$ neighbor at least one. Each point
in $Q(t-)$ has $r$ area adjacent to it, so
$b_1 + b_2 + 2 b_3 \leq\#Q(t) r$.

The variable $\pswap$ can be set to any number from $0$ to $1$: letting
$\pswap= 1/4$ gives an upper bound on the bad event rate of
$(1/2)b_1 + (1/2) b_2 + b_3 \leq(1/2)\#Q(t) r$.

Recall the bad event rate when $\pswap= 0$ was bounded above by
$\#Q(t) r$. With $\pswap= 1/4$, the bad event rate is bounded above
by $\#Q(t) r / 2$, and this factor of two carries throughout the remainder
of the proof to give (\ref{EQN:secondTheorem}).
\end{pf*}

H\"aggstr\"om and Steif gave a result similar to the previous theorem for
finitary codings for high noise Markov random
fields \cite{haggstroms2000}, but their analysis involves moving backwards
rather than forwards in time, and their result does not employ
the swap move.

\begin{figure}

\includegraphics{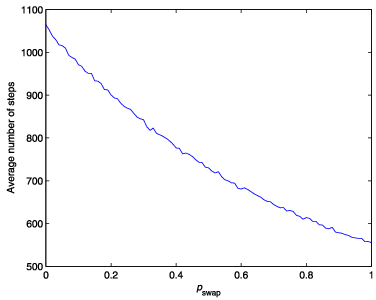}

\caption{Running time of dCFTP for Strauss model on $S = [0,1]^2$,
$\beta_1 = 50$, $\beta_2 = 0.5$, $R = 0.05$, $\lambda$ is Lebesgue measure,
as $\pswap$ runs from 0 to 1.}
\label{FIG:runtimes}
\end{figure}

Figure \ref{FIG:runtimes} illustrates the
mean run time for a fixed value of $\lambda$ as the probability of a swap
varies from $p = 0$ up to $p = 1$. The running time (as measured
by generated iterations) decreases as the chance of swapping increases.
This same phenomenon was
noted for hard core gas models on graphs \cite{huber2000a}, and at
present is unexplained by theory.

%

\section{Conclusions}
The regular birth--death chains only move when no point blocks the birth
of a point in the dominating process. The birth--death swap chains move
when at most one point blocks the birth of a point in the dominating
process. This alone means that more moves are being taken, and
helps to explain the improved analysis and improved performance
when used for perfect sampling with dominated coupling from the past.



\printhistory

\end{document}